\input amstex
\documentstyle{amsppt}
\overfullrule=0pt


\font\sc=cmcsc10

\topmatter
\title
APOLOGY OF EUCLID
\endtitle
\abstract
This is a short apology of the style of the {\it Elements\/}
by Euclid and Bourbaki.
\endabstract

\author
S.~S. Kutateladze
\endauthor
\date
April 21, 2005\enddate
\address
Sobolev Institute of Mathematics\newline
Novosibirsk State University\newline
\endaddress
\email
sskut\@member.ams.org
\endemail

\endtopmatter                       

\document
\baselineskip=.95\baselineskip
A somewhat derogatory term ``bourbakism'' proliferates in many public
discussions about the teaching of mathematics. We hear many funny
anecdotes about commutativity  as a method of calculation
as well as separate addition of nominators and denominators. 
Professional mathematicians and teachers divide into the hostile groups
that discuss with  alienation and indignation of the medieval
scholastics  the ``problem of the naturalness of zero'' 
as well as the priority rights between
the concepts ``greater than,'' ``greater than or equal to,''
and ``strictly greater than.'' All these stories and philippics  
are nice and true to some extent but  rest upon 
a clear-cut misunderstanding.

It stands to reason to recall that there was no teacher whose name was
Bourbaki. It is also reasonable to bear in mind that the treatise of
Bourbaki is written as imitation of {\it Euclid's Elements}. The style
of Bourbaki's {\it Elements of Mathematics} is exactly the style of
Euclid.

Any serious criticism of the books by Bourbaki bases on pretensions to
their content rather than style. Bourbaki's treatise is evidently
incomplete. Many important mathematical theories are absent or covered
inadequately. A few volumes present the dead ends of exuberant
theories. All these shortcomings are connected with the major capital
distinction between the books by Euclid and Bourbaki. In his {\it
Elements} Euclid set forth the theory that was almost complete in his
times, the so-called ``Euclidean'' plane and space geometry. Most of
this section of science was made clear once and forever in the epoch
of Euclid.

The Bourbaki project was implemented in the period of very rapid
progress in mathematics. Many books of the treatise became obsolete at
the exact moment of publication. In particular, functional analysis
had been developing contrary to what one might imagine reading the
book {\it Topological Vector Spaces}. But to a failure was doomed the
heroic and ambitions plan of Bourbaki to present the elements of the
whole mathematics of the twentieth century in a single treatise along
the methodological lines of Euclid. Mathematics renews and enriches
itself with outstanding brilliant achievements much faster than the
books of Bourbaki's treatise were compiled. There is no wonder
that the mathematical heroes who create  the twentieth century
mathematics have distinctly and immediately scented the shortcomings of
Bourbaki. The treatise encountered severe criticism and even
condemnation since it omits many important topics. As  usual, this
serious criticism convened all sorts of educationists,
 would-be specialists in ``propaedeutics'' and ``methodology'' who
are hardly aware of what is going on in the real mathematics.
Everyone knows that to criticize a book for incompleteness is a weak
argument since it is strange to judge an article for what is absent
in this article. Grudges against the content of the treatise transform by necessity to
the criticism of its form. The terseness, conciseness, and lapidary of
the style of exposition fall victim to criticism and even ostracism by
the adversaries of the malicious ``bourbakism'' in education.

Salomon Bochner, one of the famous mathematicians of the past, observed
with a witty smile:

\medskip \item{}{\eightpoint\sl\indent Also, if examined
 ``objectively,'' Euclid's work ought to have been any educationist's
 nightmare. The work presumes to begin
from a beginning; that is, it presupposes a certain level of
readiness, but makes no other prerequisites. Yet it never offers any
``motivations,'' it has no illuminating ``asides,'' it does not attempt
to make anything ``intuitive,'' and it avoids ``applications'' to a
fault. It is so ``humorless'' in its mathematical purism that,
although it is a book about ``Elements,'' it nevertheless does not
unbend long enough in its singlemindedness to make the remark, however
incidentally, that if a rectangle has a base of 3 inches and a height
of 4 inches then it has an area of 12 square inches. Euclid's work
never mentions the name of a person; it never makes a statement about,
or even an (intended) allusion to, genetic developments of
mathematics; it makes no cross references, except once, the exception
being in proposition 2 of Book 13, where the text refers to, and
repeats the content of, the ``first theorem of the tenth book,''
which, as it happens, is Euclid's ``substitute'' for the later axiom
of Archimedes. Euclid has a fixed pattern for the enunciation of a
proposition, and, through the whole length of 13 books, he is never
tempted to deviate from it. In short, it is almost impossible to
refute an assertion that the Elements is the work of an unsufferable
pedant and martinet... Euclid's work became one of the all-time best
sellers. According to ``objective'' Pestalozzi criteria, it should
have been spurned by students and ``progressive'' teachers in every
generation. But it nevertheless survived intact all the turmoils,
ravages, and illiteracies of the dissolving Roman Empire, of the early
Dark Ages, of the Crusades, and of the plagues and famines of the
later Middle Ages. And, since printing began, Euclid has been printed
in as many editions, and in as many languages, as perhaps no other
book outside the Bible.}\footnote {Bochner S. {\it The Role of
Mathematics in the Rise of Science.} Princeton University Press, 1981, pp.~35--36}

\medskip\noindent Euclid's book is a totally appalling, terse and
formal presentation of axioms, definitions, lemmas and theorems
without any motivation and digression, lacking any illuminating
examples from physics, economics, social or spiritual life. However,
it is the book that lives about two and a half millennia and shows no
indication of dying. In contract, the textbooks fail to survive the
gerontological tests that define the area of a figure by sowing it
with some grain or cutting it off a sheet of paper.

We must avoid mixing together the full-time and extramural forms of
training, the transfer and saving of knowledge. The Babylonian texts
on mathematics are in fact problem-books with solutions. This style of teaching is
still alive. However, no problem-book of any sort can compare with
{\it Euclid's Elements} in its long-term impact on mathematics and
culture as a whole.
Any student's notes of a mathematical course still remind us of {\it
Euclid's Elements} and its successor in style, Bourbaki's {\it Elements
of Mathematics}. 

In common parlance,   {\it bourbakism} stands for
``formalistic structural mathematics,'' whatever the bizarre term
means. In fact, this vogue word rarely implies
something more than a simple reference to the century-old tradition of
shortening and saving mathematical theories in axiomatic form. This
marvelous and noble tradition stems from the writings of Euclid. 
Elimination of extravagancy and pursuit of consistency, clarity, 
terseness,  and rationality in exposition stimulate, organize, and discipline  
mind and thought, revealing  the intrinsic beauty and harmony of mathematics.
It is exactly the impersonal style of {\it Euclid's Elements}, 
lacking any temporal inklings, that makes them especially valuable 
and allows anybody to understand what they tell us  when
centuries have elapsed.

The ``verbal'' problems, practical motivations
and emphasis on a person's creativity as well as the subjective
coloring of exposition and present-day allusions are absolutely
obligatory gadgets in the tool-kit for training. However, the
particular products of these immortal teaching tools are rather
volatile, momentary, and fragile; they often die at the spur of the
moment of their enunciation.

Science must preserve old knowledge as well as meet the challenges
of nowadays by solving the new and pending problems. Therefore,
teaching has the twofold task of preserving and transferring
knowledge, ``filling the mind'' in combination with ``igniting a
fire,'' i.e., the initiation and stimulation of creative search into
new knowledge. There is no reason to oppose the transfer and
preservation of knowledge and the training of creativity and practical
skills in raising and solving the problems of today. Preservation of
mathematical knowledge in the impersonal and dry style of textbooks
never excludes the possibility of creative search of the teacher.
On the contrary, the style of
Euclid presupposes perpetual creativity, calling the teacher for
finding and using subtle personal adjustments, subjective keys and
even mysteries for igniting students' interest in mathematics,
the understanding of its place and role in science, industry, and
other areas of public life as well as for training  skills
of  application of mathematics in practical problems.

The everlasting duty of the teacher is to destroy the obstacles to the
understanding of mathematics,  reveal the liberating essence of its
free thinking, and  explain that {\sc mathematics is the most human
of all  human sciences.} There is no math without a man or a woman.
The physical world still prevails but math vanishes without men and women.
We people do math. We do it, thinking about everyone and we do it for
everybody. The purpose and essence of mathematics reside in the freedom
it brings to us.

Mathematics welcomes everyone, combining free access, democracy, and openness
with the indisputable prohibition of any prejudice, subjectiveness,
and arbitrariness of judgements.

One of the most personalized sciences which requires everybody's
personal effort for solving a however simple arithmetical problem,
mathematics has learned  to make the complex the simple and
comprehensible to each of us.

The most human among sciences, mathematics has elaborated its
beautiful ``unhuman'' form of the objective transfer of knowledge in
writing---the classic style of the Hellenistic ``Elements.''

There are no King's ways to mathematics; the road to mathematics was
charted by Euclid. The style of Euclid not only lives in the books by
Bourbaki but also proliferates in hundreds of thousands of students'
notes throughout the world. This style is an achievement and 
article of pride of our ancient science.

\enddocument